\documentclass[12pt,a4paper]{article}
\usepackage[english]{babel}
\usepackage{amsmath}
\usepackage{amssymb}
\usepackage{graphicx}
\usepackage{color}

\usepackage{pgf}
\usepackage{pgffor}
\usepgflibrary{plothandlers}
\usepackage{pgfplots}
\pgfplotsset{compat=newest}
 \pgfplotsset{width=15cm}
\pgfplotsset{plot coordinates/math parser=false}
\newlength\figureheight
\newlength\figurewidth

\newtheorem{de}{Definition}[section]
\newtheorem{theo}{Theorem}[section]
\newtheorem{cor}[theo]{Corollary}

\newtheorem{lemm}[theo]{Lemma}
\newtheorem{rem}[theo]{Remark}
\newcommand{\Proof}{\medskip\noindent {\bf Proof. }}

\newcommand{\C}{{\cal C}}

\newcommand{\N}{\mathbb{N}}

\newcommand{\R}{\mathbb{R}}

\newcommand{\al}{\alpha}

\newcommand{\e}{\varepsilon}
\newcommand{\Om}{\Omega}

\newcommand{\s}{\sigma}

\newcommand{\h}{{H^{1}(\Om)}}

\newcommand{\hzero}{{H^1_0(\Om)}}

\newcommand{\ii}{\infty}

\newcommand{\ton}{\text{ \ on }}
\newcommand{\tin}{\text{ \ in }}

\newcommand{\nn}{\nonumber}
\newcommand{\bs}{\backslash}
\newcommand{\dr}{\partial}

\newcommand{\g}{\nabla}
\newcommand{\ti}{\tilde}
\newcommand{\p}{\cdot}
\newcommand{\ul}{\underline}
\newcommand{\ol}{\overline}
\newcommand{\su}{\subset}
\newcommand{\La}{\triangle}
\newcommand{\supp}{\text{supp}}

\newcommand{\ra}{\longrightarrow}
\newcommand{\cqfd}{\hfill $\square$\\ \medskip}
\newcommand{\n}[2]{\left\|{#1}\right\|_{#2}}

\title{A mathematical and numerical framework for ultrasonically-induced Lorentz force electrical impedance tomography\thanks{\footnotesize This work was supported  by the
ERC Advanced
Grant Project MULTIMOD--267184.}}
\author{Habib Ammari\thanks{\footnotesize Department of Mathematics and Applications,
Ecole Normale Sup\'erieure, 45 Rue d'Ulm, 75005 Paris, France
(habib.ammari@ens.fr, pierre.millien@ens.fr,
laurent.seppecher@ens.fr).} \and Pol Grasland-Mongrain\thanks{\footnotesize
Inserm, U1032, LabTau, universit\'e de Lyon,  Lyon, F-69003, France
(pol.grasland-mongrain@ens-cachan.fr).}
\and Pierre Millien\footnotemark[2] \and Laurent Seppecher\footnotemark[2] \and
Jin-Keun Seo\thanks{\footnotesize Department of Computational
Science and Engineering, Yonsei University, 50 Yonsei-Ro,
Seodaemun-Gu, Seoul 120-749, Korea (seoj@yonsei.ac.kr).}}

\begin{document}
\maketitle
\begin{abstract}
We provide a mathematical analysis and a numerical framework for Lorentz force electrical conductivity imaging. Ultrasonic vibration of a tissue in the presence of a static magnetic field induces an electrical current by the Lorentz force. This current can be detected by electrodes placed around the tissue; it is proportional to the velocity of the ultrasonic pulse, but depends nonlinearly on the conductivity distribution.
The imaging problem is to reconstruct the conductivity distribution from measurements of the induced current. To solve this nonlinear inverse problem, we first make use of a virtual potential to relate explicitly the current measurements to the conductivity distribution and the velocity of the ultrasonic pulse. Then, by applying a Wiener filter to the measured data, we
reduce the problem to imaging the conductivity from an internal electric current density. We first introduce an optimal control method for solving such a problem.
A new direct reconstruction scheme involving a partial differential equation is then proposed based on viscosity-type regularization to a transport equation satisfied by the current density field. We prove that solving such an equation yields the true conductivity distribution as the regularization parameter approaches zero.  We also test both schemes numerically in the presence of measurement noise, quantify their stability and resolution, and compare their performance.
\end{abstract}

\bigskip

\noindent {\footnotesize Mathematics Subject Classification
(MSC2000): 35R30, 35B30.}


%%%%%%%%%%%%%%%%%%%%%%%%%%%%%%%%%%%%%%%%%%%%%%%%%%%%%%%%%%%%%%%%%%%%%%%%%%%%%
\section{Introduction}
%%%%%%%%%%%%%%%%%%%%%%%%%%%%%%%%%%%%%%%%%%%%%%%%%%%%%%%%%%%%%%%%%%%%%%%%%%%%%
Ultrasonic imaging is currently used in a wide range of medical diagnostic applications. Its high spatial resolution, combined with a real-time imaging capability, lack of side effects, and relatively low cost make it an attractive technique. However, it can be difficult to differentiate soft tissues because acoustic impedance varies by less than $10\%$ among muscle, fat, and blood \cite{goss1978comprehensive}. In contrast, electrical conductivity varies widely among soft tissue types and pathological states \cite{Foster1989, morimoto1993study} and its measurement can provide
information about the physiological and pathological condition of tissue \cite{laure}.
Several techniques have been developed to map electrical conductivity. The most well known is electrical impedance tomography, in which electrodes are placed around the organ of interest, a voltage difference is applied, and the conductivity distribution can be reconstructed from the measurement of the induced current at the electrodes \cite{ammari2008introduction,tscan, cheney1999electrical}. This technique is  harmless to the patient if low currents are used. However, the ill-posed character of the inverse problem results in lower spatial resolution than that achieved by ultrasound imaging, and any speckle information is lost.

The Lorentz force plays a key role in acousto-magnetic tomographic techniques \cite{roth}.
Several approaches have been developed with the aim of providing electrical impedance information at a spatial resolution on the scale of ultrasound wavelengths  \cite{ejam, GraslandMongrain2013, he1, he2, montalibet2002etude, roth, roth2009ultrasonically, wen1998hall}. These include Hall effect imaging, magneto-acoustic current imaging, magneto-acoustic tomography with magnetic induction, and ultrasonically-induced Lorentz force imaging. Acousto-magnetic tomographic techniques have the potential to detect small conductivity inhomogeneities, enabling them to diagnose pathologies such as cancer by detecting tumorous tissues when other conductivity imaging techniques fail to do so.

\begin{figure}[!h]
\def\svgwidth{\linewidth}
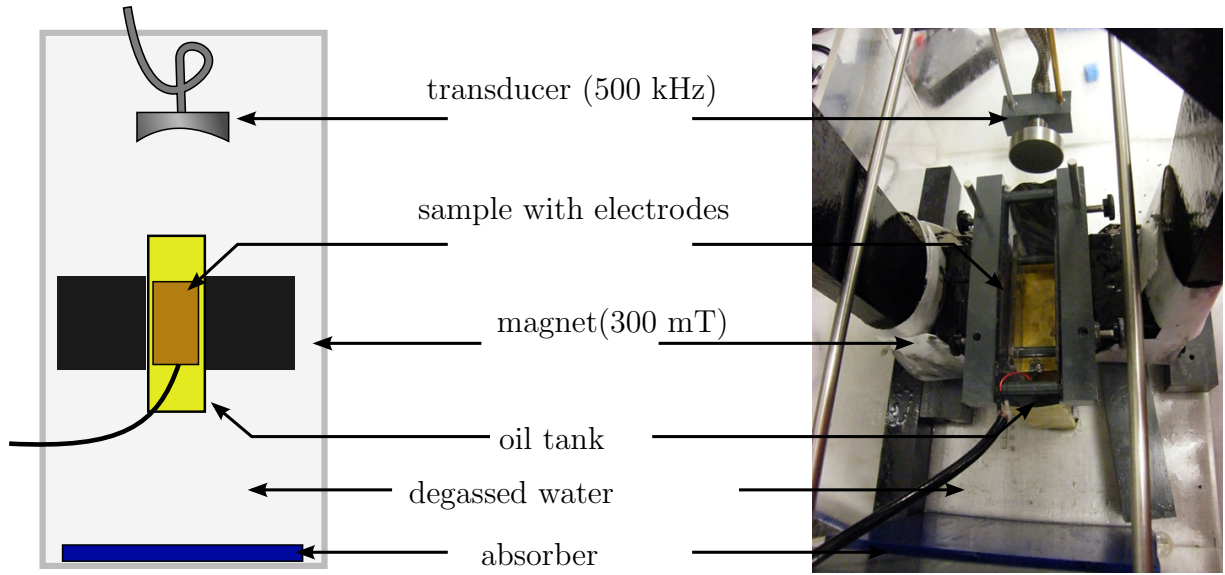
\caption{Example of the imaging device. A transducer is emitting ultrasound in a sample placed in a constant magnetic field. The induced electrical current is collected by two electrodes.}\label{pictureexperiment}
\end{figure}

In ultrasonically-induced Lorentz force method (experimental apparatus presented in Figure~\ref{pictureexperiment}) an ultrasound pulse propagates through the medium to be imaged in the presence of a static magnetic field. The ultrasonic wave induces Lorentz' force on the ions in the medium, causing the negatively and positively charged ions to separate. This separation of charges acts as a source of electrical current and potential. Measurements of the induced current give information on the conductivity in the medium. A $1$ \emph{Tesla} magnetic field and a $1$ \emph{MPa} ultrasonic pulse
induce current at the \emph{nanoampere} scale. Stronger magnetic fields and ultrasonic beams can be used to enhance the signal-to-noise ratio \cite{GraslandMongrain2013}.

This paper provides a rigorous mathematical and numerical framework for
ultra\-sonically-induced Lorentz force electrical impedance tomography.
We develop two efficient methods for reconstructing the conductivity in the medium from the induced electrical current.
As far as we know, this is the first mathematical and numerical modeling of the experiment conducted in \cite{GraslandMongrain2013} to illustrate the feasibility of ultrasonically-induced Lorentz force electrical impedance tomography. Earlier attempts to model mathematically this technique were made in \cite{ejam, ip}.

The paper is organized as follows. We start by describing the ionic model of conductivity. From this model we derive the current density induced by an ultrasonic pulse in the presence of a static magnetic field. We then find an expression of the measured current. The inverse problem is to image the conductivity distribution from such measurements corresponding to different pulse sources and directions.
A virtual potential used with simple integrations by parts can relate the measured current to the conductivity distribution and the velocity of the ultrasonic pulse. A Wiener deconvolution filter can then reduce the problem to imaging the conductivity from the internal electric current density.
The internal electric current density corresponds to that which would be induced by a constant voltage difference between one electrode and another with zero potential.
We introduce two reconstruction schemes for solving the imaging problem from the internal data. The first is an optimal control method; we also propose an alternative to this scheme via the use of a transport equation satisfied by the internal current density. The second algorithm is direct and can be viewed as a PDE-based reconstruction scheme. We prove that solving such a PDE yields to the true conductivity distribution as the regularization parameter tends to zero.
In doing so,  we prove the existence of the characteristic lines for the transport equation under some conditions on the conductivity distribution.
We finally test numerically the two proposed schemes in the presence of measurement noise, and also quantify their stability and resolution.

The ultrasonically-induced Lorentz force electrical impedance tomography investigated here can be viewed as a new hybrid technique for conductivity imaging. It has been experimentally tested \cite{GraslandMongrain2013}, and was reported to produce images of quality comparable to those of ultrasound images taken under similar conditions. Other emerging hybrid techniques for conductivity imaging have also been reported \cite{eric, ejam, ip12, yves, otmar1, sirev, otmar2}.

%%%%%%%%%%%%%%%%%%%%%%%%%%%%%%%%%%%%%%%%%%%%%%%%%%%%%%%%%%%%%%%%%%%%%%%%%%%%%
\section{Electric measurements from acousto-magnetic coupling}
%%%%%%%%%%%%%%%%%%%%%%%%%%%%%%%%%%%%%%%%%%%%%%%%%%%%%%%%%%%%%%%%%%%%%%%%%%%%%
Let a physical object to be imaged occupy a three-dimensional domain $\Om$ with a smooth boundary $\dr\Om$. Assume that this body is placed in a constant magnetic field $B$ in the direction $e_3$ where $\{e_1,e_2,e_3\}$ denotes the standard orthonormal basis of $\R^3$. We are interested in recovering the electrical conductivity of this body $\s\in L^\ii(\Om)$ with the known lower and upper bounds:
\begin{equation}\nn
0<\ul\s\leq\s\leq\ol \s<\ii.
\end{equation}
An acoustic transducer sends a short acoustic pulse from $y\in\R^3$ in the direction $\xi\in S^2$, with $S^2$ being the unit sphere, such that $\xi\cdot e_3=0.$
This pulse generates the velocity field $v(x,t)\xi$ with $v(x,t)$ taking the following form:
\begin{equation}\label{eq_pulse}\begin{aligned}
v(x,t)=w\big(z-ct\big)~A\big(z,|r|\big),
 \end{aligned}\end{equation}
where
$$z=(x-y)\cdot\xi\quad\mbox{and }\quad r=x-y-z\xi~\in \Upsilon_{\xi}:= \{ \zeta\in \R^3~:~\zeta\cdot \xi=0\}.
$$
Here, $w\in \C^\ii_c\big(\R\big)$, supported in $]-\eta,0[$, is the ultrasonic pulse profile;  $A\in \C^\ii\big(\R\times\R^+\big)$, supported in $\R^+\times[0,R]$, is the cylindrical profile distribution of the wave corresponding to the focus of the acoustic transducer; and $R$ is the maximal radius of the acoustic beam.

%%%%%%%%%%%%%%%%%%%%%%%%%%%%%%%%%%%%%%%%%%%%%%%%%%%%%%%%%%%%%%%%%%%%%%%%%%%%%%%%%%%%%%%%%%
\subsection{The ionic model of conductivity}

We describe here the electrical behavior of the medium as an electrolytic tissue composed of ions capable of motion in an aqueous tissue. We consider $k$ types of ions in the medium with charges of $q_i$, $i\in\{1,\dots,k\}$. The corresponding volumetric density $n_i$ is assumed to be constant. Neutrality in the medium is described as

\begin{equation}\label{chargen} \begin{aligned}
\sum_{i}q_in_i=0.
\end{aligned}\end{equation}

The Kohlrausch law defines the conductivity of such a medium as a linear combination of the ionic concentrations

\begin{equation}\label{eq_sigma}
\sigma=e^+\sum_i \mu_i q_i n_i,
\end{equation}
where $e^+$ is the elementary charge, and the coefficients $\mu_i$ denote the ionic mobility of each ion $i$. See, for example, \cite{montalibet2002etude,pride1994governing}.

%%%%%%%%%%%%%%%%%%%%%%%%%%%%%%%%%%%%%%%%%%%%%%%%%%%%%%%%%%%%%%%%%%%%%%%%%%%%%%%%%%%%%%%%%%
\subsection{Ion deviation by Lorentz force}

We embed the medium in a constant magnetic field $B$ with direction $e_3$, and perturb it mechanically using the short, focused, ultrasonic pulses $v$ defined in (\ref{eq_pulse}). The motion of the charged particle $i$ inside the medium is deviated by the Lorentz force
\begin{equation}
F_i=q_iv \xi\times B.
\end{equation}
This force accelerates the ion in the orthogonal direction $\tau=\xi\times e_3$. Then, almost immediately, the ion reaches a constant speed given by

\begin{equation}\nn
v_{\tau,i}= \mu_i |B|v
\end{equation}
at the first order. See \cite{montalibet2002etude, pride1994governing} for more details. Finally, the ion $i$ has a total velocity
\begin{equation}\nn\begin{aligned}
v_i =v\xi+\mu_i  |B| v  \tau.\\
\end{aligned}\end{equation}
The current density generated by the displacement of charges can be described as follows:
\begin{equation}\nn\begin{aligned}
j_S =\sum_i n_iq_iv_i=\left(\sum_i n_i q_i\right)v\xi+\left(\sum_i n_i
\mu_i q_i \right)|B|v \tau.\\
\end{aligned}\end{equation}
Using the neutrality condition (\ref{chargen}) and the definition of $\s$ in (\ref{eq_sigma}), we get the following simple formula for $j_S$:

\begin{equation}\label{def_js}\begin{aligned}
j_S=\frac{1}{e^+}|B|\sigma v\tau ,
\end{aligned}\end{equation}
which is in accordance with the formula used in \cite{ejam}.

This electrolytic description of the tissue characterizes the interaction between the ultrasonic pulse and the magnetic field through a small deviation of the charged particles embedded in the tissue. This deviation generates a current density $j_S$ orthogonal to $\xi$ and to $B$, locally supported inside the domain. At a fixed time $t$, $j_S$ is supported in the support of $x\mapsto v(x,t)$. This current is proportional to $\sigma$, and is the source of the current that we measure on the electrodes placed at $\dr\Om$. In the next section, a formal link is substantiated between $j_S$ and the measured current $I$.

%%%%%%%%%%%%%%%%%%%%%%%%%%%%%%%%%%%%%%%%%%%%%%%%%%%%%%%%%%%%%%%%%%%%%%%%%%%%%%%%%%%%%%%%%%
\subsection{Internal electrical potential}

\begin{figure}[!h]
\def\svgwidth{\linewidth}
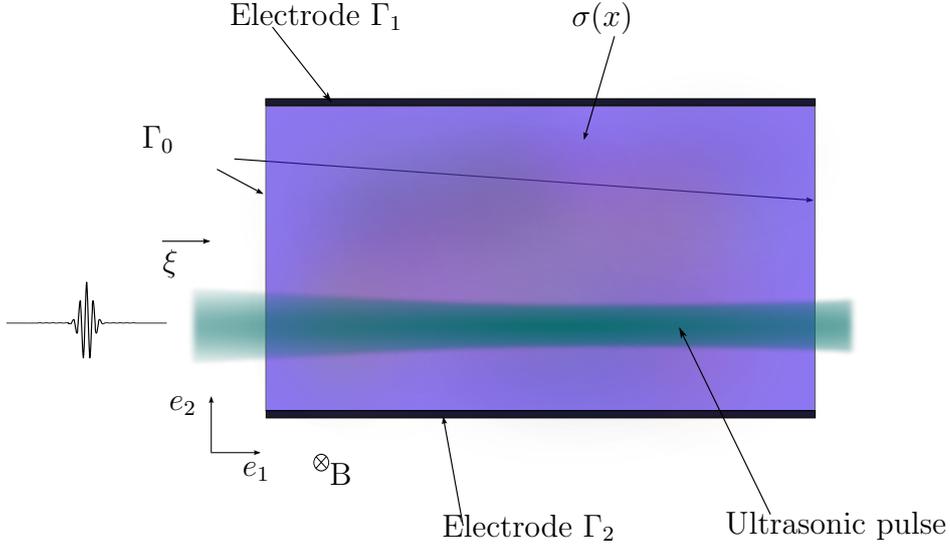
\caption{Imaging system configuration. An ultrasonic wave propagates in a medium of electrical conductivity $\sigma$ comprised between electrodes $\Gamma_1$ and $\Gamma_2$.}\label{exemple_config}
\end{figure}

Because the characteristic time of the acoustic propagation is very long compared with the electromagnetic wave propagation characteristic time, we can adopt the electrostatic frame.
Consequently, the total current $j$ in $\Omega$ at a fixed time $t$ can be formulated as
\begin{equation}\begin{aligned}
j =j_S+\sigma\g u,
\end{aligned}\end{equation}
where $u$ is the electrical potential. It satisfies
\begin{equation}\begin{aligned}
\g\p(j_S+\sigma\g u)=\g\p j =0.
\end{aligned}\end{equation}
Figure \ref{exemple_config} shows the configuration under consideration. Let $\Gamma_1$ and $\Gamma_2$ be portions of the boundary $\partial \Om$ where two planner electrodes are placed. Denote $\Gamma_0=\partial\Om\setminus (\Gamma_1\cup\Gamma_2)$.

As we measure the current between the two electrodes $\Gamma_1$ and $\Gamma_2$, the electrical potential is the same on both electrodes, and can be fixed to zero without loss of generality.
Further, it is assumed that no current can leave from $\Gamma_0$. The potential $u$ can then be defined as the unique solution in $\h$ of the elliptic system
\begin{equation}\label{eq_cond}\left\{\begin{aligned}
-\g\p(\s\g u)&=\g\p j_S&\tin \Om,\\
u&=0&\ton \Gamma_1\cup\Gamma_2,\\
\dr_\nu u&=0&\ton \Gamma_0.
\end{aligned}\right.\end{equation}
Throughout this paper $\dr_\nu$ denotes the normal derivative. Note that the source term $j_S$ depends on the time $t>0$, the position of the acoustic transducer $y\in\R^3$, and the direction $\xi\in S^2$. The electrical potential $u$ also depends on these variables. 

The measurable intensity $I$ is the current flow through the electrodes. Integrating (\ref{eq_cond}) by parts gives
\begin{equation}\nn\begin{aligned}
\int_{\Gamma_1}\s\dr_\nu u+\int_{\Gamma_2}\s\dr_\nu u=0,
\end{aligned}\end{equation}
which is the expression of current flow conservation. We define the intensity $I$ by
\begin{equation}\label{eq_I}\begin{aligned}
I=\int_{\Gamma_2}\s\dr_\nu u.
\end{aligned}\end{equation}

%%%%%%%%%%%%%%%%%%%%%%%%%%%%%%%%%%%%%%%%%%%%%%%%%%%%%%%%%%%%%%%%%%%%%%%%%%%%%%%%%%%%%%
\subsection{Virtual potential}
In order to link  $I$  to $\sigma$, we introduce a virtual potential
$U\in\h$  defined as the unique solution of
\begin{equation}\label{eq_U}
\left\{\begin{aligned}
-\g\p(\s\g U) &= 0 &\tin\Om,\\
U&=0 &\ton\Gamma_1,\\
U &=1 &\ton\Gamma_2,\\
\dr_\nu U&=0&\ton\Gamma_0.
\end{aligned}\right.
\end{equation}
Then we multiply (\ref{eq_cond}) by $U$ and integrate by parts. Assuming that the support of $v$ does not intersect the electrodes $\Gamma_1$ and $\Gamma_2$, we obtain
\begin{equation}\nn\begin{aligned}
-\int_\Om\s\g u\p\g U+\int_{\Gamma_2}\s\dr_\nu u=\int_\Om j_S\p\g U.
\end{aligned}\end{equation}
From the property of $U$ in  (\ref{eq_U}) and the definition of $I$  in  (\ref{eq_I}), the above identity becomes
\begin{equation}\nn\begin{aligned}
I=\int_\Om j_S\p\g U.
\end{aligned}
\end{equation}
The above identity links the measured intensity $I$ to an internal information of $\sigma$ using the expression of $j_S$ in (\ref{def_js}):
\begin{equation}\nn\begin{aligned}
I=\frac{|B|}{e^+}\int_\Om  v(x,t) \sigma(x)\g U(x)dx \p\tau.
\end{aligned}
\end{equation}
According to (\ref{eq_pulse}), $v$ depends on $y$, $\xi$, and $t$, so does $I$. We define the measurement function as
\begin{equation}\begin{aligned}
M_{y,\xi}(z) =\int_\Om v(x,z/c)\s(x)\g U(x)dx\p\tau(\xi)
\end{aligned}
\end{equation}
for any $y\in\R^3$, $\xi\in S^2$ and $z>0$. We assume the knowledge of this function in a certain subset of $\R^3\times S^2\times\R^+$ denoted by $Y\times \mathfrak{S}\times ]0,z_{max}[$. We will discuss later the assumptions  we have to impose on this subset in order to make the reconstruction accurate and stable.

%%%%%%%%%%%%%%%%%%%%%%%%%%%%%%%%%%%%%%%%%%%%%%%%%%%%%%%%%%%%%%%%%%%%%%%%%%%%%%%%%%%%%%%%%%%
\section{Construction of the virtual current}
%%%%%%%%%%%%%%%%%%%%%%%%%%%%%%%%%%%%%%%%%%%%%%%%%%%%%%%%%%%%%%%%%%%%%%%%%%%%%%%%%%%%%%%%%%%

For simplicity, let us restrict ourselves to the two dimensional case where both the conductivity $\sigma$ and the virtual potential $U$ do not change in $e_3$-direction. For convenience, the same notations will be used as in the three dimensional case.

In order to obtain the information of $\sigma$ contained in $M_{y,\xi}$, we need to separate the contribution of the displacement term $v$ from this measurement function. Using the cylindrical symmetry of this integration we write for any $z\in]0,z_{max}[$,

\begin{equation}\label{msigdec} \begin{aligned}
M_{y,\xi}(z)&=\int_\R\int_{\Upsilon_\xi}  w(z-z')(\s\g U)(y+z'\xi+r)A(z',|r|)drdz'\p\tau(\xi),\\
&=\int_\R w(z-z')\int_{\Upsilon_\xi} (\s\g U)(y+z'\xi+r)A(z',|r|)dr dz'\p\tau(\xi), \\
&=\left(W\star\Phi_{y,\xi}\right)(z)\p\tau(\xi),
\end{aligned}\end{equation}
where $W(z)=w(-z)$,  $\star$ denotes the convolution product, and

\begin{equation}\nn
\Phi_{y,\xi}(z)=\int_{\Upsilon_\xi} \s(y+z\xi+r)A(z,|r|)\g U(y+z\xi+r) dr.
\end{equation}

As will be shown in section \ref{numerical}, through a one dimensional deconvolution problem that can be stably solved using, for instance, a Wiener-type filtering method, we get access to the function $\Phi_{y,\xi}\p\tau(\xi)$. Now the question is about the reconstruction of $\s$ from $\Phi_{y,\xi}\p\tau(\xi)$. We can notice that $\Phi_{y,\xi}$ is a weighted Radon transform applied to the virtual current field $\s\g U$. The weight $A(z,|r|)$ is critical for the choice of the method that we can use. Closer this weight is to a Dirac mass function, better is the stability of the reconstruction. In this case, if the field $\s\g U$ does not have too large variations, we can recover a first-order approximation; as discussed in the rest of this section.

In order to make the reconstruction accurate and stable, we make two assumptions on the set of parameters $Y\times D\times]0,z_{max}[$. For any $x\in\Om$, we define
\begin{equation}\nn
\mathfrak{S}_x=\left\{\xi\in \mathfrak{S}~:~\xi=\frac{x-y}{|x-y|}~\mbox{ for some }~y\in Y\right\}.
\end{equation}
The first assumption is

\begin{equation}\nn
\mbox{(H1)} ~~~\ \forall x\in\Om,\ \ \exists\ \xi_1,\xi_2 \in \mathfrak{S}_x ~~\mbox{ s.t. }~~ |\xi_1\times\xi_2|\not=0,
\end{equation}
 and the second one reads
\begin{equation}\nn
\mbox{(H2)}~~~ \forall x\in\Om,\ \ \forall \xi\in \mathfrak{S}_x,\ \ \exists \mbox{ unique}~~ y\in Y ~~\mbox{ s.t. }~~ \xi=\frac{x-y}{|x-y|}.
\end{equation}

From the assumption (H2), we can define a distance map $|x-y|$ as a function of $x$ and $\xi$. We will denote $d_Y(x,\xi)=|x-y|$. By a change of variables, we rename our data function $\Sigma$ as

\begin{equation}\label{leading1} \begin{aligned}
\psi(x,\xi) &=\Phi_{y,\xi}\big(d_Y(x,\xi)\big)\p\tau(\xi)\\
&=\int_{\Upsilon_\xi} (\s\g U)(x+r)A\big(d_Y(x,\xi),|r|\big) dr\p\tau(\xi).
\end{aligned}\end{equation}
Now if we denote by \begin{equation}
\label{leading2}
\gamma(x,\xi) = \int_{\Upsilon_\xi}A\big(d_Y(x,\xi),|r|\big) dr\ \tau(\xi),
\end{equation} then we expect that
$$\psi(x,\xi)\approx (\s\g U)(x)\p \gamma(x,\xi),$$
provided the $\supp(A)$ is small enough and $\s\g U$ does not vary too much.  The following lemma makes this statement precise.

\begin{lemm} \label{lembv} Consider a fixed direction $\xi\in \mathfrak{S}$ and consider the domain covered by the pulses of direction $\xi$ defined by $\Om_\xi=\{x\in\Om~:~ \xi\in \mathfrak{S}_x\}$. Suppose that the virtual current $\s\g U$ has bounded variations, then
\begin{equation}\nn
\n{\psi(\cdot,\xi)-\s\g U\cdot\gamma(\cdot,\xi)}{L^1(\Om_\xi)}\leq cR \|\s\g U\|_{TV(\Om)^2},
\end{equation}
where $R$ is the maximum radius of the cylindrical support of the envelope $A$ and $c>0$ depends on the shape of $A$. Here, $\|\;\|_{TV(\Om)^2}$ denotes the total variation semi-norm.
\end{lemm}
\Proof For a.e. $x\in \Om_\xi$, we have
\begin{equation}\nn\begin{aligned}
&\left|\psi(x,\xi)-(\s\g U)(x)\cdot\gamma(x,\xi)\right| \leq\\
&\int_{\Upsilon_\xi} \left|(\s\g U)(x+r)-(\s\g U)(x)\right|A\big(d_Y(x,\xi),|r|\big) dr,
\end{aligned}\end{equation}
and so
\begin{equation}\nn\begin{aligned}
&\n{\psi(\cdot,\xi)-\s\g U\cdot\gamma(\cdot,\xi)}{L^1(\Om_\xi)}\\
&\leq \int_{\Upsilon_\xi} \int_{\Om_\xi}\left|(\s\g U)(x+r)-(\s\g U)(x)\right|A\big(d_Y(x,\xi),|r|\big) dxdr\\
&\leq \|\s\g U\|_{TV(\Om)^2}\int_{\Upsilon_\xi} |r|\sup_{0<z<z_{max}}A(z,|r|) dr\\
&\leq  2\pi R \|\s\g U\|_{TV(\Om)^2}\int_{\R+} \sup_{0<z<z_{max}}A(z,\rho) d\rho.\\
\end{aligned}\end{equation}\cqfd

Note that in the most interesting cases, $\s\g U$ has bounded variations.
For example, if $\s$ has a piecewise $W^{1,\ii}$ smoothness on smooth inclusions, then $\s\g U$ has bounded variations. This also holds true for $\sigma$ in some subclasses  of functions of bounded variations. In the following, we make the assumption, as in Lemma \ref{lembv}, that  $\s\g U$ has bounded variations.

In conclusion, our data approximates the quantity $(\s\g U)(x)\cdot\gamma(x,\xi)$ for any $x\in\Om$, $\xi\in \mathfrak{S}_x$ where the vector $\gamma(x,\xi)$ is supposed to be known. To get the current $(\s\g U)(x)$, we simply consider data from two linearly independent directions. Using assumption (H1), for a fixed $x\in\Om$, there exist $\xi_1,\xi_2 \in \mathfrak{S}_x$ such that $\det(\xi_1,\xi_2)\not= 0$. We construct the  $2\times 2$ invertible matrix

\begin{equation}\nn
\Gamma(x,\xi_1,\xi_2)=\left[\begin{matrix}
\gamma(x,\xi_1)^\perp\\
\gamma(x,\xi_2)^\perp
\end{matrix}\right],
\end{equation}
and the data column vector

\begin{equation}\nn
\Psi(x,\xi_1,\xi_2)=\left[\begin{matrix}
\psi(x,\xi_1)\\
\psi(x,\xi_2)
\end{matrix}\right].
\end{equation}

We approximate the current $\sigma\g U(x)$ by the vector field

\begin{equation}\nn
V(x,\xi_1,\xi_2)=\Gamma(x,\xi_1,\xi_2)^{-1} \Psi(x,\xi_1,\xi_2).
\end{equation}
Indeed, for any open set $\widetilde\Om\su\Om_{\xi_1}\cap\Om_{\xi_2}$, the following estimate holds:

\begin{equation}\nn\begin{aligned}
&\n{V(\cdot,\xi_1,\xi_2)-\s\g U}{L^1(\widetilde\Om)^2} \\
&\leq \sup_{x\in \widetilde\Om}\n{\Gamma(x,\xi_1,\xi_2)^{-1}}{{\cal L}(\R^2)}\left(\sum_{i=1}^2\n{\psi(\cdot,\xi_i)-\s\g U\p\gamma(\cdot,\xi_i)}{L^1(\Om_{\xi_i})}\right)^{1/2}\\
&\leq cR \|\s\g U\|_{TV(\Om)^2}.
\end{aligned}\end{equation}

It is worth mentioning that if more directions are available, then we can use them to enhance the stability of the reconstruction. The linear system becomes over-determined and we can get the optimal approximation by using a least-squares method.

%%%%%%%%%%%%%%%%%%%%%%%%%%%%%%%%%%%%%%%%%%%%%%%%%%%%%%%%%%%%%%%%%%%%%%%%%%
\section{Recovering the conductivity by optimal control}
%%%%%%%%%%%%%%%%%%%%%%%%%%%%%%%%%%%%%%%%%%%%%%%%%%%%%%%%%%%%%%%%%%%%%%%%%%

In this section we assume that, according to the previous one, we are in the situation where we know a good approximation of the virtual current $D:=\s\g U$ in the sense of $L^1(\Om)^2$. The objective here is to provide efficient methods for separating $\s$ from $D$.

For $a<b$, let us denote by $L^\ii_{a,b}(\Om):=\{f\in L^\ii(\Om)~:~\ a<f<b\}$ and define the operator $\mathcal F:L^\ii_{\ul\s ,\ol\s}(\Om)\ra\h$ by

\begin{equation}\label{eq_U2}
\mathcal F[\s]=U:\left\{\begin{aligned}
\g\p(\s\g U) &= 0 &\tin\Om,\\
U&=0 &\ton\Gamma_1,\\
U &=1 &\ton\Gamma_2,\\
\dr_\nu U&=0&\ton\Gamma_0.
\end{aligned}\right.
\end{equation}
The following lemma holds.
\begin{lemm} The operator $\mathcal F$ is Fr\'echet differentiable and for any $\s\in L^\ii_{\ul\s ,\ol\s}(\Om)$ and $h\in L^\ii(\Om)$ such that $\s+h\in L^\ii_{\ul\s ,\ol\s}(\Om)$ we have

\begin{equation}\label{eqdeffh}
d\mathcal F[\s](h)=v:~~\left\{\begin{aligned}
\g\p(\s\g v) &=-\g\p (h\g \mathcal F[\s])  &\tin\Om,\\
v&=0 &\ton\Gamma_1\cup\Gamma_2,\\
\dr_\nu v&=0&\ton\Gamma_0.
\end{aligned}\right.
\end{equation}
\end{lemm}

\Proof Let us denote by $w=\mathcal F[\s+h]-\mathcal F[\s]-v$. This function is in $\h$ and satisfies the equation

\begin{equation}\nn
\g\p(\s\g w) =-\g\p (h\g (\mathcal F[\s+h]-\mathcal F[\s]))
\end{equation}
with the same boundary conditions as $v$. We have the elliptic global control:
\begin{equation}\nn
\|\g w\|_{L^2(\Om)}\leq\frac{1}{\ul\s}\n h {L^\ii(\Om)}\n{\g (\mathcal F[\s+h]-\mathcal F[\s])}{L^2(\Om)}.
\end{equation}
Since
\begin{equation}\nn
\g\p(\s\g (\mathcal F[\s+h]-\mathcal F[\s])) =-\g\p (h\g \mathcal F[\s +h]),
\end{equation}
we can also control $\mathcal F[\s+h]-\mathcal F[\s]$ with
\begin{equation}\nn
\n{\g (\mathcal F[\s+h]-\mathcal F[\s])}{L^2(\Om)}\leq\frac{1}{\sqrt{\ul\s}}\n h {L^\ii(\Om)}\n{\g \mathcal F[\s+h]}{L^2(\Om)}.
\end{equation}
Then, there is a positive constant $C$ depending only on $\Om$ such that
 $$\n{\g \mathcal F[\s+h]}{L^2(\Om)}\leq  C \sqrt{\frac{\ol\s}{\ul\s}}.$$
Finally, we obtain
\begin{equation}\nn
\n {\g w}{L^2(\Om)}\leq C\frac{\sqrt{\ol\s}}{\ul\s^2}\n h {L^\ii(\Om)}^2.
\end{equation}\cqfd

We look  for the minimizer of the functional

\begin{equation}
J[\s]=\frac 1 2\int_\Om\left|\s \g \mathcal F[\s]-D\right|^2.
\end{equation}
In order to do so, we compute its gradient. The following lemma holds.

\begin{lemm} \label{lemcompj} For any $\s\in L^\ii_{\ul\s ,\ol\s}(\Om)$,

\begin{equation}\nn
d J[\s]=(\s\g \mathcal F[\s]-D-\g p)\p \g \mathcal F[\s],
\end{equation}
where $p$ is defined as the solution to the adjoint problem:

\begin{equation} \label{defp}
\left\{\begin{aligned}
\g\p(\s\g p) &=\g\p (\s^2\g \mathcal F[\s]-\s D)  &\tin\Om,\\
p&=0 &\ton\Gamma_1\cup\Gamma_2,\\
\dr_\nu p&=0&\ton\Gamma_0.
\end{aligned}\right.
\end{equation}
\end{lemm}

\Proof As $\mathcal F$ is Fr\'echet differentiable, so is $J$. For $\s\in L^\ii_{\ul\s ,\ol\s}(\Om)$ and $h\in L^\ii(\Om)$ such that $\s+h\in L^\ii_{\ul\s ,\ol\s}(\Om)$, we have
\begin{equation}\nn
d J[\s](h)=\int_\Om(\s \g \mathcal F[\s]-D)\p(h \g \mathcal F[\s]+\s \g d\mathcal F[\s](h)).
\end{equation}
Now, multiplying (\ref{defp}) by $d\mathcal F[\s](h)$, we get

\begin{equation}\nn
\int_\Om\s\g p\p\g d\mathcal F[\s](h)=\int_\Om(\s^2\g \mathcal F[\s]-\s D)\p\g d\mathcal F[\s](h).
\end{equation}
On the other hand, multiplying (\ref{eqdeffh}) by $p$ we arrive at
\begin{equation}\nn
\int_\Om\s\g p\p\g d\mathcal F[\s](h)=-\int_\Om h\g \mathcal F[\s]\p\g p,
\end{equation}
and therefore,
\begin{equation}\nn
d J[\s](h)=\int_\Om h(\s \g \mathcal F[\s]-D-\g p)\p\g \mathcal F[\s].
\end{equation}\cqfd

Lemma \ref{lemcompj} allows us to implement a numerical gradient descent method in order to find $\s$. A regularization term can also be added to $J[\s]$ in order to avoid instability. As we are seeking discontinuous $\s$ with smooth variations out of the discontinuity set, a good choice would be the minimization of the regularized functional:

\begin{equation}
J_\e[\s]=\frac 1 2\int_\Om\left|\s \g \mathcal F[\s]-D\right|^2+\e ||\s ||_{TV(\Om)},
\end{equation}
where $\e >0$ is the regularization parameter.

%See, for instance, \red{give references}.

%%%%%%%%%%%%%%%%%%%%%%%%%%%%%%%%%%%%%%%%%%%%%%%%%%%%%%%%%%%%%%%%%%%%%%%%%%
\section{The orthogonal field method}
In this section, we present an alternative direct method to optimal control for reconstructing the conductivity $\sigma$ from the internal data $\sigma \nabla U$. It is based on solving a transport equation.
The following approach may be extended to the three dimensional case. However,  several proofs would need to be revisited.

Given a vector field $D=\s\g U$ which is parallel to $\g U$ everywhere, we may construct the vectorial field $F=(D_2,-D_1)$ which is everywhere orthogonal to $D$. The flow of $F$ may define the level sets of $U$. Assuming that the variations of the conductivity $\s$ are far enough from $\Gamma_0$, we can assume that $U(x)=x_2$ on this boundary part. Then $U$ is a solution of the following transport equation:

\begin{equation}\label{eq_trans}\left\{\begin{aligned}
F\p\g u &=0\ &\tin\Om,\\
u &=x_2\ &\ton\dr\Om.
\end{aligned}\right.
\end{equation}
In the case where (\ref{eq_trans}) is well posed and can be solved, we can reconstruct the virtual potential $U$. The conductivity $\sigma$ is deduced from $U$ and $D$ by the following identity
\begin{equation}\label{defsigma}\begin{aligned}
\sigma=\frac{D \cdot \g U}{\vert D\vert^2}.
\end{aligned}\end{equation}
Despite to its very simple form, this first-order equation is really tricky. Existence and uniqueness are both difficult challenges in the general case. Our main difficulty here is due to the fact that $F$ is discontinuous. As the function $U$ that we are looking for is a natural solution of this equation, we are only concerned here with the uniqueness of a solution to (\ref{eq_trans}).

%%%%%%%%%%%%%%%%%%%%%%%%%%%%%%%%%%%%%%%%%%%%%%%%%%%%%%%%%%%%%%%%%%%%%%%%%%%%%%%%%%%%%%%%%%%%%%%%%%%%%%%%%%%%%%%%%%%%%%%%%%%%%%%%%%%%%%%
\subsection{Uniqueness result for the transport equation}

The uniqueness of a solution to (\ref{eq_trans}) is directly linked to the existence of outgoing characteristic lines defined by the dynamic system:

\begin{equation}\label{char}\left\{\begin{aligned}
X'(t) &=F(X(t)),\ t\geq 0,\\
X(0) &=x, \ x \in\Om,
\end{aligned}\right.\end{equation}
which usually needs the continuity of $F$. As $\s$ is in general not continuous,  $F$ is not continuous, which makes the classical existence results useless. Nevertheless, under some assumptions on $\sigma$, we can insure the existence of the characteristic lines.

\begin{de} For any $k\in\N$, $\al\in ]0,1[$, for any curve $\mathcal C$ of class $C^{1,\alpha}$ such that  $\Om\setminus \mathcal C$ is a union of connected domains $\Om_i, i=1,2,\cdots n$,   we define $C_{\mathcal C}^{k,\al}\big(\ol\Om\big)$ to be the class of functions $f:\Om\ra\R$ satisfying

\begin{equation}\nn\begin{aligned}
 f|_{\Om_i}\in C^{k,\al}\big(\ol{\Om_i}\big)\quad \forall i=1,\cdots n.
\end{aligned}\end{equation}
\end{de}

\begin{de} A conductivity $\sigma$ is said to be admissible if there exists a constant $\al\in ]0,1[$ and a curve $\mathcal C$ of class $C^{1,\al}$ such that
$\s\in C_{\mathcal C}^{0,\al}\big(\ol\Om\big)\cap L^\ii_{\ul\s ,\ol\s}(\Om)$ and
\begin{equation}\nn\begin{aligned}
\inf_{\Om\setminus\mathcal C} \s\g \mathcal F[\sigma]\p e_2~ >~0.
\end{aligned}\end{equation}
\end{de}
 If $\s$ is admissible and belongs to $C_{\mathcal C}^{0,\al}\big(\ol\Om\big)$, then the solution $U$ of (\ref{eq_U}) belongs to $C_{\mathcal C}^{1,\al}\big(\ol\Om\big)$ and the field $F=(\s\g U)^\perp$ satisfies
\begin{equation}\nn\begin{aligned}
&F\in C_{\mathcal C}^{0,\al}\big(\ol\Om\big)\quad\mbox{and}\quad \inf_{\Om\setminus\mathcal C} F\cdot e_1~>~0.
\end{aligned}\end{equation}

Moreover, as $F$ is orthogonal to $\sigma\g U$, we can describe the jump of $F$  at the curve $\mathcal C$. Defining the normal and tangential unit vectors $\nu$ and $\tau$ and also the local sides (+) and (-) with respect to $\nu$, we can write $F$ on both sides as

\begin{equation}\nn\begin{aligned}
F^+=\s^+\dr_\nu U^+\tau+\s^+\dr_\tau U^+\nu,\\
F^-=\s^-\dr_\nu U^-\tau+\s^-\dr_\tau U^-\nu\\
\end{aligned}\end{equation}
with the transmission conditions, $\s^+\dr_\nu U^+=\s^-\dr_\nu U^-$ and $\dr_\tau U^+=\dr_\tau U^-$. Finally, we characterize the discontinuity of $F$ by

\begin{equation}\nn\begin{aligned}
 \left[F\right]=[\s]\dr_\tau U\nu,
\end{aligned}\end{equation}
where $[\, ]$ denotes the jump across $\mathcal C$.

With all of these properties for the field $F$, we can prove the existence of the characteristic lines for (\ref{char}).

\begin{theo}\label{theo_local}(Local existence of characteristics) Assume that $F\in C_{\mathcal C}^{0,\alpha}\big(\ol\Om\big)$ with $\mathcal C$ of class $C^{1,\alpha}$ for $\alpha\in ]0,1[$. Assume that the discontinuity of $F$ on $\mathcal C$ satisfies
\begin{equation}\nn\begin{aligned}
F^+=f\tau+\s^+g\nu,\\
F^-=f\tau+\s^-g\nu
\end{aligned}\end{equation}
with $f,g,\s^+,\s^-\in C^{0,\alpha}(\mathcal C)$ where $\s^+,\s^-$ are positive and $g$ is locally signed.
Then, for any $x_0\in\Om$, there exists $T>0$ and $X\in C^1\big([0,T[,\Om\big)$ such that $t\mapsto F(X(t))$ is measurable and

 \begin{equation}\nn\begin{aligned}
X(t)=x_0+\int_0^tF(X(s))ds,\ \ \forall t\in[0,T[.
\end{aligned}\end{equation}
\end{theo}

\Proof If $x_0\notin \mathcal C$, then $F$ is continuous in a neighborhood of $x_0$ and  the Cauchy-Peano theorem can be applied.

If $x_0\in \mathcal C$, then we choose a disk $B\subset \Om$ centered at $x_0$. The oriented line $\mathcal C$ separates $B$ in two simply connected open domains called $B^+$ and $B^-$.
For ease of explanation, we may assume that $\mathcal C\cap B$ is straight line (since we can flatten the curve using a proper $C^{0,\alpha}$-diffeomorphism).

Assume that $g(x_0)>0$. Up to rescaling $B$, we can assume that $g(x)>0$ for all $x \in \mathcal C \cap B$.   We extend $F|_{B^+}$ to a continuous field $\ti{F}\in C^{0}(B)$ by even reflection.
The Cauchy-Peano theorem insures the existence of $T>0$ and $X\in C^1\big([0,T[,\Om\big)$ such that $X(0)=x_0$ and $X'(t)=\ti{F}(X(t))$ for all $t\in[0,T[$. 
As $g(x_0)>0$, we have $X'(0)\p\nu(x_0)>0$ and $X(t)\in \ol{B^+}$ in a neighborhood of $0$. Thus, for a small enough $t$, $X'(t)=F(X(t))$. If $g(x_0)<0$, then we apply the same argument by interchanging $B^-$ and $B^+$.

Suppose now that $g(x_0)=0$. The field $F$ is now tangent to the discontinuity line.
If $f(x_0)=0$, then $X(t)=x_0$ is a solution. We assume here that $f(x_0)>0$. 
As $g$ is assumed to be locally signed, we can suppose that $g\geq 0$ in a small sub-curve of $\mathcal C$ satisfying $(x-x_0)\p\tau(x_0)>0$.  Again, we extend $F|_{B^+}$ to a continuous field $\ti{F}\in C^{0}(B)$ by even reflection and use the Cauchy-Peano theorem to show that there exists $T>0$ and $X\in C^1\big([0,T[,\Om\big)$ such that $X(0)=x_0$ and $X'(t)=\ti {F}(X(t))$ for all $t\in [0,T[$. In order to complete the proof, we should show that $X(t)$ belongs to $\ol{B^+}$ for $t$ small enough. If not, there exists a sequence $t_n\searrow 0$ such that $X(t_n)\in B^-$. By the mean value theorem, there exists $\tilde t_n\in (0, t_n)$ such that $F(X(\tilde t_n))\cdot\nu(x_0)=X'(\tilde t_n)\cdot\nu(x_0)<0$.   Thus, $X(t)$ belongs to $\ol{B^+}$ and $X'(t)=F(X(t))$ for $t$ small enough.

Note that the local monotony of $g$ is satisfied in many cases. For instance if $\mathcal C$ is analytic and $\sigma$ is piecewise constant, then $\nabla U$ is analytic on $\mathcal C$ and hence, $g$ is locally signed. \cqfd

It is worth mentioning that existence of a solution for the Cauchy problem
(\ref{char}) has been proved in \cite{bressan} provided that $F\cdot \nu >0$  on $\mathcal C$. Here, we have made a weaker assumption. In fact, we only need that $F\cdot \nu$ is locally signed.

\begin{cor}\label{cor_maximal}(Existence of outgoing characteristics) Consider $F\in C_{\mathcal C}^{0,\alpha}(\Om)$ satisfying the same conditions as in Theorem \ref{theo_local} and the condition

\begin{equation}\nn\begin{aligned}
\inf_{\Om\setminus\mathcal C}F\p e_1\geq c,
\end{aligned}\end{equation}
where $c$ is a positive constant. Then for any $x_0\in\Om$ there exists $0<T<T_{\max}$ where $T_{\max}=\dfrac 1 c\text{diam}(\Om)$ and $X\in C^0\big([0,T[,\Om\big)$ satisfying

\begin{equation}\nn\begin{aligned}
&X(t)=x_0+\int_0^t F(X(s))ds,\ \ \forall t\in [0,T[,\\
&\lim_{t\to T}X(t)\in\dr\Om.
\end{aligned}\end{equation}
This result means that from any point $x_0\in\Om$, the characteristic line reaches $\dr\Om$ in a finite time.
\end{cor}

\Proof Let $x_0\in\Om$ and $X\in\C^0\big([0,T[,\Om\big)$ a maximal solution of (\ref{char}). Using $F\p e_1\geq c$ we have that $X'(t)\p e_1\geq c$ and so $X(t)\p e_1\geq x_0\p e_1+ct$ and as $X(t)\in\Om$ for all $t\in[0,T[$, it is necessary that $T<T_{\max}$. As $F\in C^{0,\alpha}_{\mathcal C}(\Om)$, $F$ is bounded, $X$ is Lipschitz, and the limit of $X(t)$ when $t$ goes to $T$ exists in $\ol\Om$ and is called $X(T)$. Let us show that $X(T)\in\dr\Om$. Suppose that $X(T)\in\Om$, then applying Theorem \ref{theo_local} at $X(T)$, we can continuously extend $X$ on $[T,T+\e[$ for some positive $\e$ which contradicts the fact that $X$ is a maximal solution.\cqfd

\begin{cor}\label{cor_unique}(Uniqueness for the transport problem) Consider $F\in C_{\mathcal C}^{0,\alpha}(\Om)$ satisfying the same conditions as in Corollary \ref{cor_maximal} and consider $u\in C^0\big(\ol\Om\big)\cap C^1_\mathcal C\big(\ol\Om\big)$. If $u$ is a solution of the system

\begin{equation}\label{eq_trans0}\left\{\begin{aligned}
F\p\g u &=0\ &\tin\Om,\\
u &=0\ &\ton\dr\Om,
\end{aligned}\right.\end{equation}
then $u=0$ in $\Om$.
\end{cor}

\Proof Consider $x_0\in\Om$ and a characteristic $X\in\C^0\big([0,T[,\Om\big)$ satisfying

\begin{equation}\nn\begin{aligned}
&X(t)=x_0+\int_0^t F(X(s))ds,\ \ \forall t\in [0,T[,\\
&\lim_{t\to T}X(t)\in\dr\Om.
\end{aligned}\end{equation}
We define $f\in C^0\big([0,T],\R\big)$ by $f(t)=u(X(t))$. We show that $f$ is constant. Let us define $I=X^{-1}(\mathcal C)$ then $f$ is differentiable in $[0,T]\bs I$ and $f'(t)=\g u(X(t))\p F(X(t))=0$. Let us take $t\in I$. If $t$ is not isolated in $I$, using the fact that $\partial_\tau u^+$ and $\partial_\tau u^-$
are locally signed,
$F(X(t))$ is parallel to $\mathcal C$ and for an $\e>0$, $X(s)\in \ol{B^+}$ (or $\ol{B^-})$ for $s\in[t,t+\e[$. Then, $f(s)=u(x(s))$ is differentiable on $[t,t+\e[$ with $f'(s)=\g u^+(X(s))\p F(X(s))$. This proves that $f$ is right differentiable at $t$ and $(f')^+(t)=0$. By the same argument, $f$ is left differentiable at $t$ and $(f')^-(t)=0$ and so $f$ is differentiable at $t$ with $f'(t)=0$. Finally, except for a zero measure set of isolated points, $f$ is differentiable on $[0,T]$ and $f'=0$ almost everywhere. This is not enough to conclude because there exists continuous increasing functions whose derivative is zero almost everywhere. Since for all $t,s\in[0,T]$, $$|f(t)-f(s)|\leq\sup|\g u||X(t)-X(s)|\leq\sup|\g U|\sup |F||t-s|,$$
$f$ is Lipschitz and thus absolutely continuous which implies, since $f'=0$ a.e., that $f$ is constant on $[0,T]$. We finally have $u(x_0)=f(0)=f(T)=u(X(T))=0$.\cqfd

Hence we conclude that if $\sigma$ is admissible, then $U$ is the unique solution to (\ref{eq_trans}) and we can recover $\sigma$ by (\ref{defsigma}).

\begin{rem} The characteristic method can be used to solve the transport problem. However, it suffers from poor numerical stability which is exponentially growing with the distance to the boundary. To avoid this delicate numerical issue, we propose a regularized approach for solving (\ref{eq_trans}). Our approach consists in forming from (\ref{eq_trans}) a second-order PDE and adding to this PDE a small elliptic term of order two.
\end{rem}

\subsection{The viscosity-type regularization}

In this subsection we introduce a viscosity approximation to (\ref{eq_trans}).
Let $\e>0$. We regularize the transport equation (\ref{eq_trans}) by considering the well-posed elliptic problem

\begin{equation}\label{eq_visco}\left\{\begin{aligned}
\g\p \left[\left(\e I+FF^T\right)\g u_\e\right] &=0\ &\tin\Om,\\
u_\e &=x_2\ &\ton\dr\Om.
\end{aligned}\right.
\end{equation}
The main question is to understand the behavior of $u_\e$ when $\e$ goes to zero. Or more precisely,  whether $u_\e$ converges to the solution $U$ of the transport equation (\ref{eq_trans}) for a certain topology. The following result holds. 

\begin{theo} \label{thmcv}
The sequence $(u_\e-U)_{\e>0}$ converges strongly to zero in $\hzero$.
\end{theo}

\Proof We first prove that the sequence $(u_\e-U)_{\e>0}$ converges weakly to zero in $\hzero$  when $\e$ goes to zero.
 For any $\e >0$, $\ti u_\e := u_\e -U$ is in $\hzero$ and satisfies

\begin{equation}\label{viscopreuve}
\g\p \left[\left(\e I+FF^T\right)\g \ti u_\e\right] = -\e \La U \ \tin\Om\\.
\end{equation}
Multiplying this equation by $\ti u_\e $ and integrating by parts over $\Om$, we obtain

\begin{equation}\label{eq_viscoenergy}\begin{aligned}
\e\int_\Om|\g\ti u_\e|^2+\int_\Om|F\p\g\ti u_\e|^2=-\e\int_\Om \g U\p\g \ti u_\e
\end{aligned}\end{equation}
and so,

\begin{equation}\nn\begin{aligned}
\n{\ti u_\e}{\hzero}^2\leq\int_\Om|\g u\p\g \ti u_\e|\leq \n{U}\h\n{\ti u_\e}{\hzero}.
\end{aligned}\end{equation}
Then $\n{\ti u_\e}{\hzero}\leq \n{U}\h$. The sequence $(u_\e)_{\e>0}$ is bounded in $\hzero$ and so by  Banach-Alaoglu's theorem, we can extract a subsequence which converges weakly to $u^*$ in $\hzero$.
Multiplying (\ref{viscopreuve}) by $u^*$ and integrating by parts, we get

\begin{equation}\nn\begin{aligned}
\int_\Om \left(F\cdot \g \ti u_\e\right)\left(F\cdot \g  u^*\right)=-\e \int_\Om \g U \cdot \g u^*  - \e \int_\Om \g \ti u_\e \cdot \g u^*.
\end{aligned}\end{equation}
Taking the limit when $\e$ goes to zero,
\begin{equation}\nn\begin{aligned}
\n{F\cdot\g u^*}{L^2(\Om)} =0.
\end{aligned}\end{equation}
So $u^*$ is a solution of the transport equation (\ref{eq_trans0}), and by Corollary \ref{cor_unique}, $u^*=0$ in $\Om$. Actually, there is no need for the extraction of a subsequence to get the weak convergence result. Indeed, zero is the only accumulation point for $u_\e$ for the weak topology.
Consider a subsequence $u_{\phi(\e)}$. It is still bounded in $\hzero$. Therefore, using the same argument as above, zero is an accumulation point of this subsequence.

Now, we are ready to prove the strong convergence.  From (\ref{eq_viscoenergy}) we get that

\begin{equation}\nn\begin{aligned}
\int_\Om|\g\ti u_\e|^2\leq-\int_\Om \g U\p\g \ti u_\e,
\end{aligned}\end{equation}
and as $\ti u_\e\rightharpoonup 0$ in $\hzero$, the term in the right-hand side goes to zero when $\e$ goes to zero. Hence, $\n{\ti u_\e}{\hzero}\rightarrow 0$.\cqfd

Finally, using Theorem \ref{thmcv}  we define the approximate resistivity by

\begin{equation}\nn\begin{aligned}
\frac 1{\s_\e}=\frac{D\p\g u_\e}{|D|^2},
\end{aligned}\end{equation}
which strongly converges to $\dfrac 1 \s$ in $L^2(\Om)$.

\section{Numerical results}\label{numerical}

In this section we first discuss the deconvolution step. Then we test both the optimal control and the orthogonal  field reconstruction schemes.

\subsection{Deconvolution}
In this subsection, we consider the problem of recovering $\Phi_{y,\xi}$ from the measurements
$ M_{y,\xi}$ in the presence of noise. From (\ref{msigdec}), it is easy to see that this can be done by deconvolution. However,  deconvolution is a numerically  very unstable process.
In order to render stability we use a Wiener filter \cite{stephane}. We assume that the signal $M_{y,\xi}(.)$ is perturbed by a random white noise:
\begin{equation}\label{defmtilde}
\widetilde{M}_{y,\xi}(z)= M_{y,\xi}(z) + \mu (z),
\end{equation} where $\mu$ is a white Gaussian noise with variance $\nu^2$. Equation (\ref{defmtilde}) can be written as
\begin{equation}\nn
\widetilde{M}_{y,\xi}(z)=\left(W \star \Phi_{y,\xi}\right)(z) + \mu (z).
\end{equation}
Denote by $S(\Sigma) = \int_{\mathbb{R}} \vert \mathcal{F}({\Phi_{y,\xi}})(\omega)\vert d\omega$ the mean spectral density of $\Sigma$, where $\mathcal{F}$ is the Fourier transform.
The Wiener deconvolution filter can be written in the frequency domain as
\begin{equation}\nn
\widehat{L}(\omega)= \frac{\overline{\mathcal{F}({W})}(\omega) }{\vert \mathcal{F}({W})\vert^2(\omega) +\frac{\nu}{S(\Sigma)}}.
\end{equation} The quotient ${\nu}/{S(\Sigma)}$ is the signal-to-noise ratio. So, in order to use the filter, we need to have an a priori estimate of the signal-to-noise ratio. We then recover $\Sigma$ up to a small error by
\begin{equation}\nn
\widetilde{\Sigma}_{y,\xi}=\mathcal{F}^{-1}\left( \mathcal{F}(\widetilde{M}) \widehat{L}\right).
\end{equation}
\begin{figure}
\begin{center}
\input{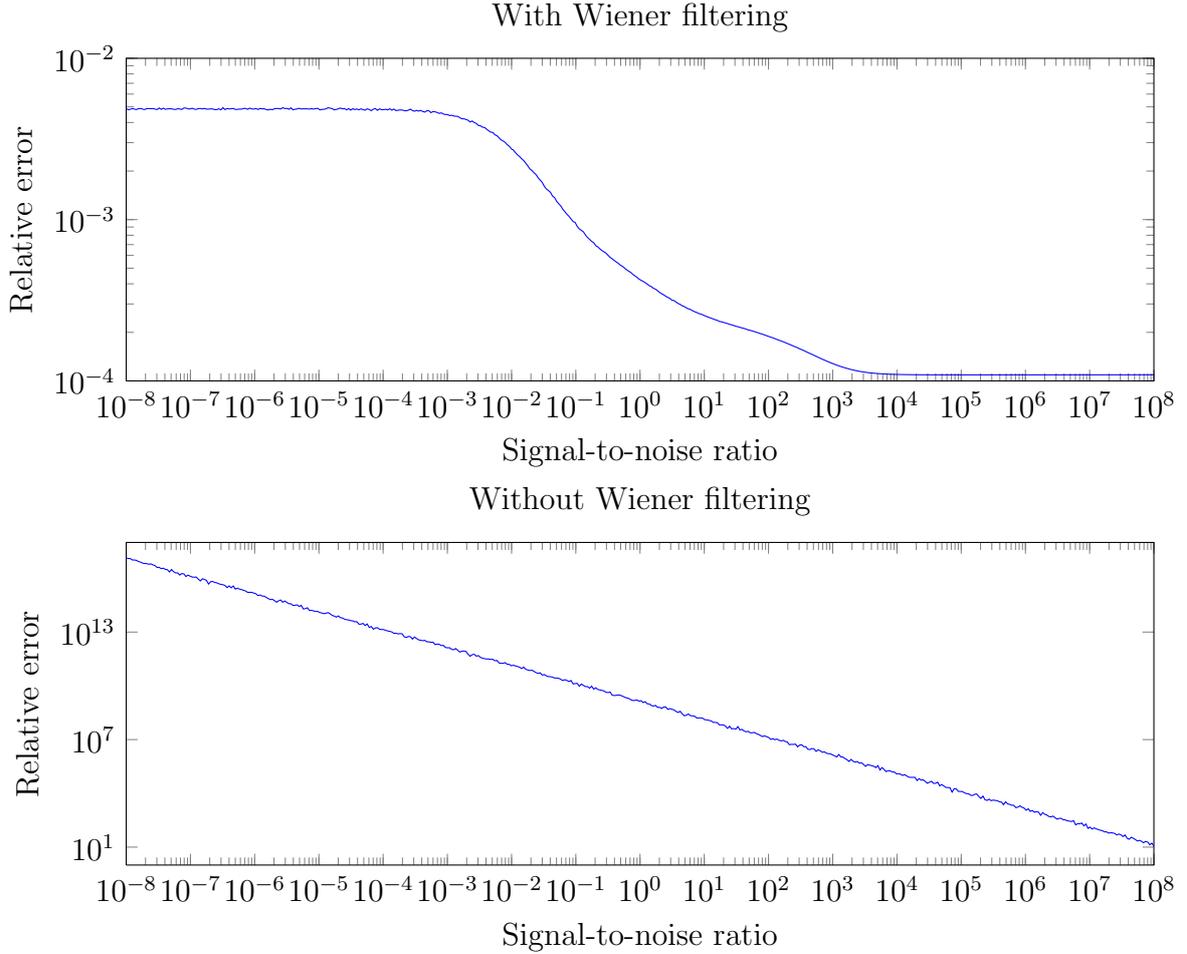}
\caption{\label{grapherror} $L^2$ norm of the relative error $\frac{\Vert \Sigma - \widetilde{\Sigma}\Vert_2}{\Vert \Sigma \Vert_2}$ with respect to the signal-to-noise ratio.}
\end{center}
\end{figure}

\subsection{Conductivity reconstructions}
\begin{figure}
\begin{center}
% This file was created by matlab2tikz v0.4.2.
% Copyright (c) 2008--2013, Nico Schlömer <nico.schloemer@gmail.com>
% All rights reserved.
% 
% The latest updates can be retrieved from
%   http://www.mathworks.com/matlabcentral/fileexchange/22022-matlab2tikz
% where you can also make suggestions and rate matlab2tikz.
% 
% 
% 
\begin{tikzpicture}

\begin{axis}[%
width=5.30555555555556in,
height=2.65277777777778in,
axis on top,
scale only axis,
xmin=0,
xmax=2,
ymin=0,
ymax=1,
colormap/jet,
colorbar,
point meta min=1,
point meta max=8
]
\addplot graphics [xmin=-0.000977517106549365,xmax=2.00097751710655,ymin=-0.000978473581213307,ymax=1.00097847358121] {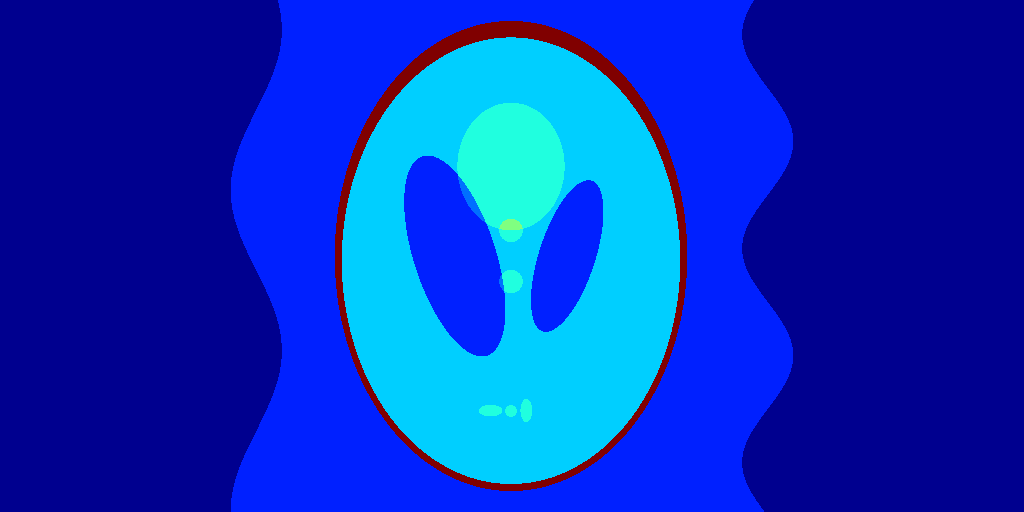};
\end{axis}
\end{tikzpicture}%
\caption{\label{conductivity} Conductivity map to be reconstructed.}
\end{center}
\end{figure}

\begin{figure}
\begin{center}
\includegraphics[scale=0.25]{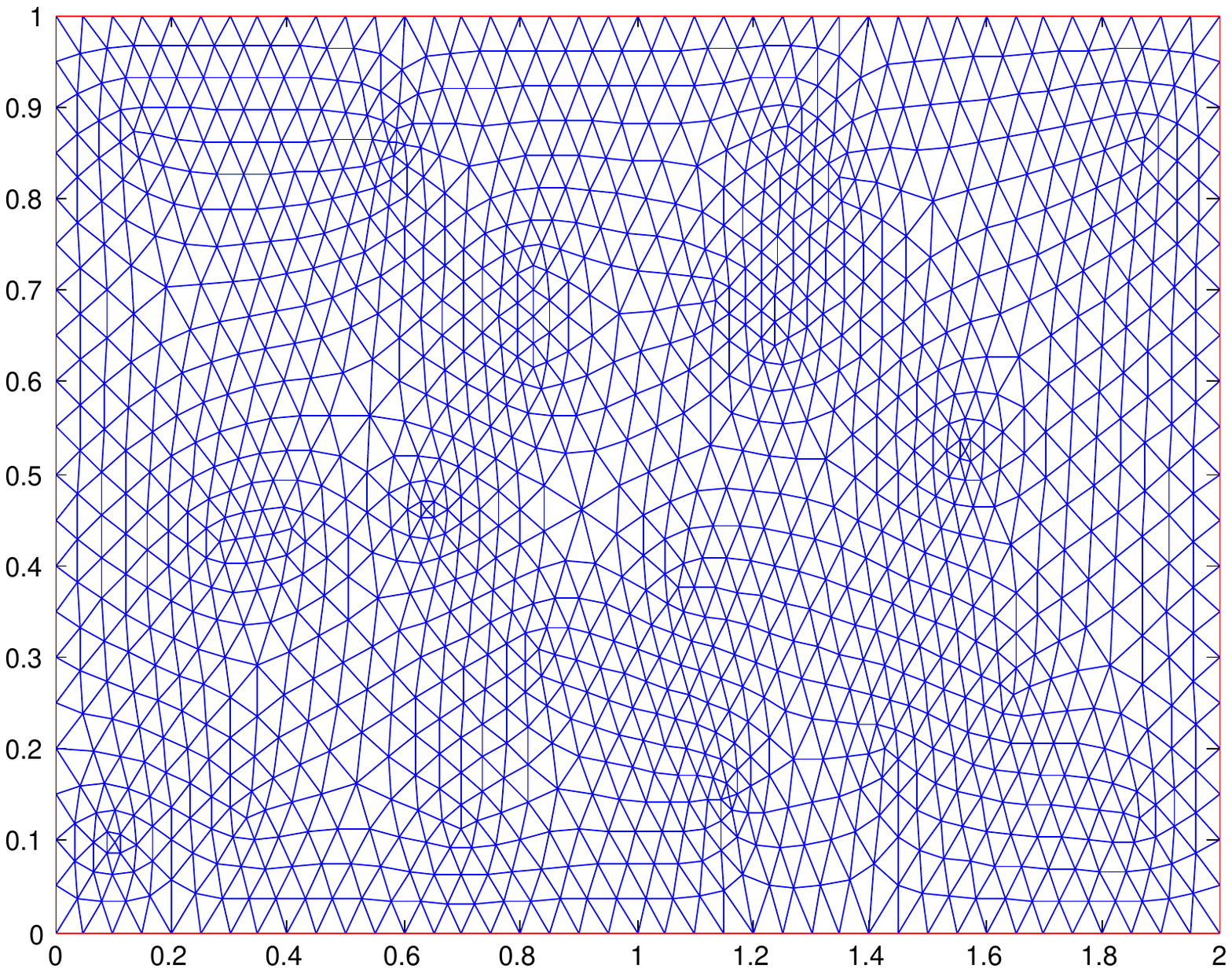}
\includegraphics[scale=0.25]{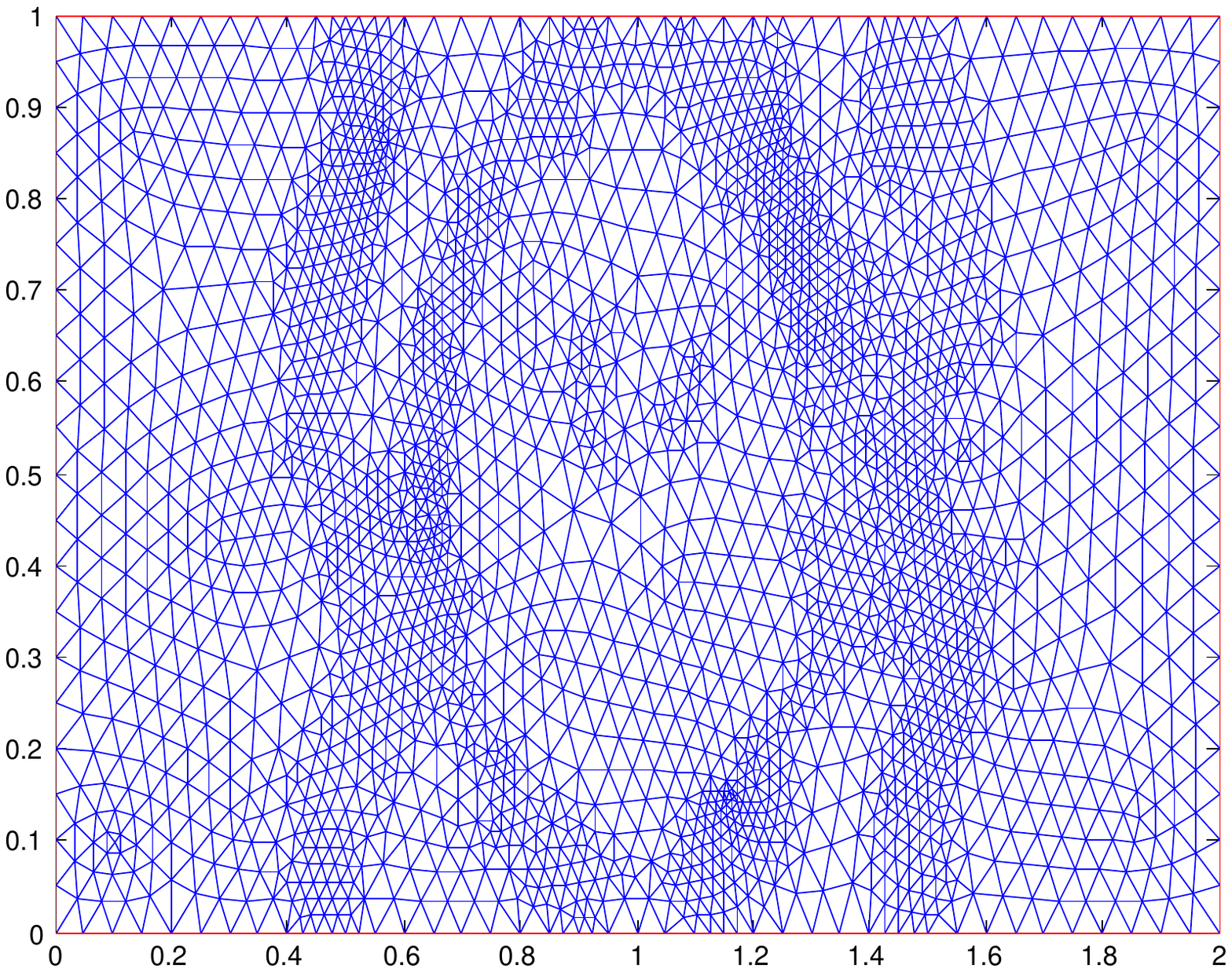}
\includegraphics[scale=0.25]{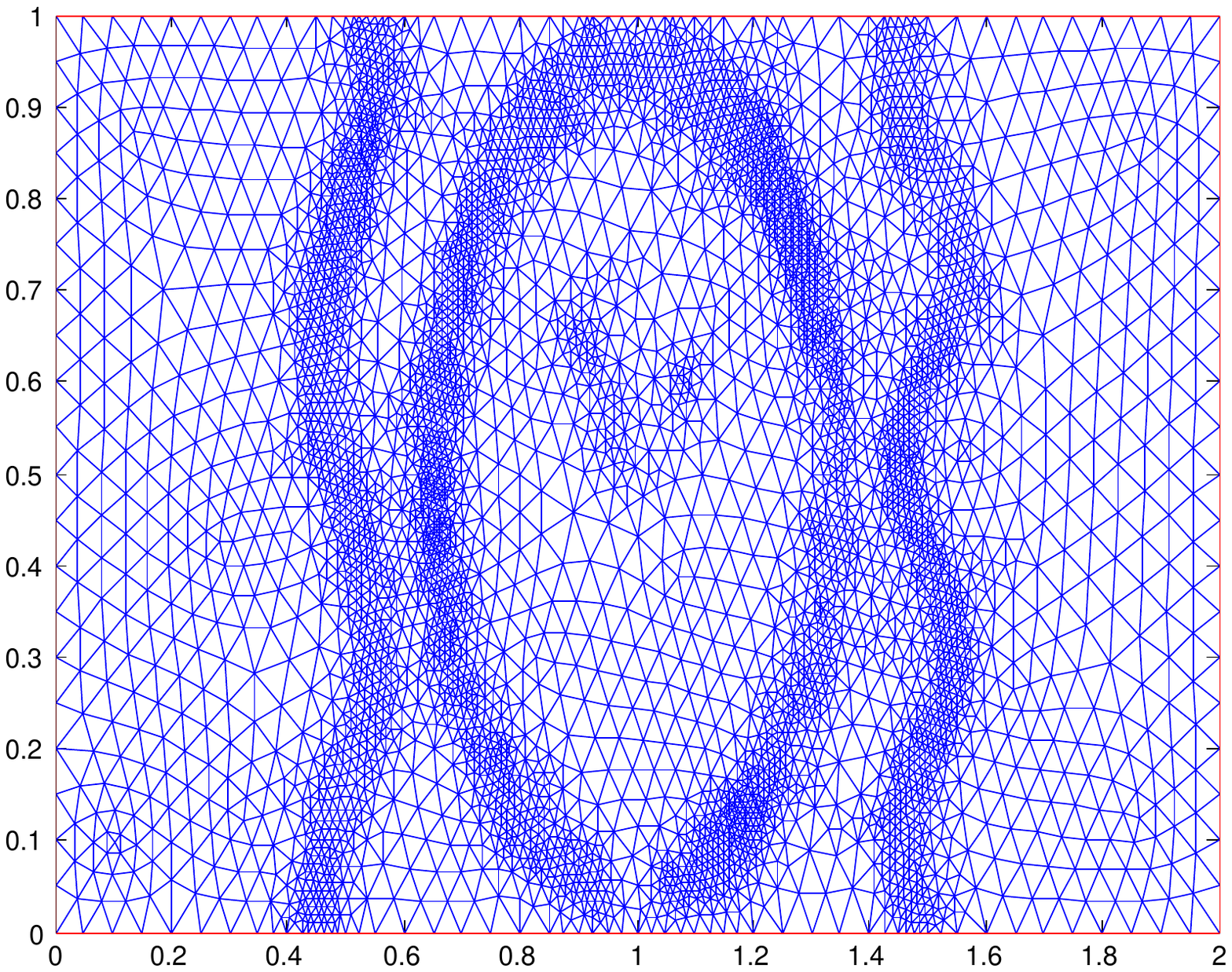}
\caption{ Meshes for solving (\ref{eq_visco}): initial mesh (left), adapted mesh (middle), and the mesh used to generated the data (right).\label{figmesh}}
\end{center}
\end{figure}

In the numerical simulations, we choose $\Omega=]0,2[\times ]0,1[$.  Figure~\ref{conductivity} shows the true conductivity map in the medium. The simulations are done using a PDE solver.
The data is simulated numerically on a fine mesh. For the orthogonal field method, in order to solve (\ref{eq_visco}), we use a coarse mesh. Then we reconstruct an initial image of the conductivity. Based on the initial image, an adaptive mesh refinement for solving (\ref{eq_visco}) yields a conductivity image of a better quality. Figure \ref{figmesh} shows the used meshes for solving the viscosity approximation.

\subsubsection{The optimal control method}
\begin{figure}
\begin{center}
% This file was created by matlab2tikz v0.4.2.
% Copyright (c) 2008--2013, Nico Schlömer <nico.schloemer@gmail.com>
% All rights reserved.
% 
% The latest updates can be retrieved from
%   http://www.mathworks.com/matlabcentral/fileexchange/22022-matlab2tikz
% where you can also make suggestions and rate matlab2tikz.
% 
% 
% 
\begin{tikzpicture}

\begin{axis}[%
width=5.30555555555556in,
height=2.65277777777778in,
axis on top,
scale only axis,
xmin=0,
xmax=2,
ymin=0,
ymax=1,
title={Optimal control method},
colormap/jet,
colorbar,
point meta min=1.48194083548821,
point meta max=2.64701330597566
]
\addplot graphics [xmin=-0.000977517106549365,xmax=2.00097751710655,ymin=-0.000978473581213307,ymax=1.00097847358121] {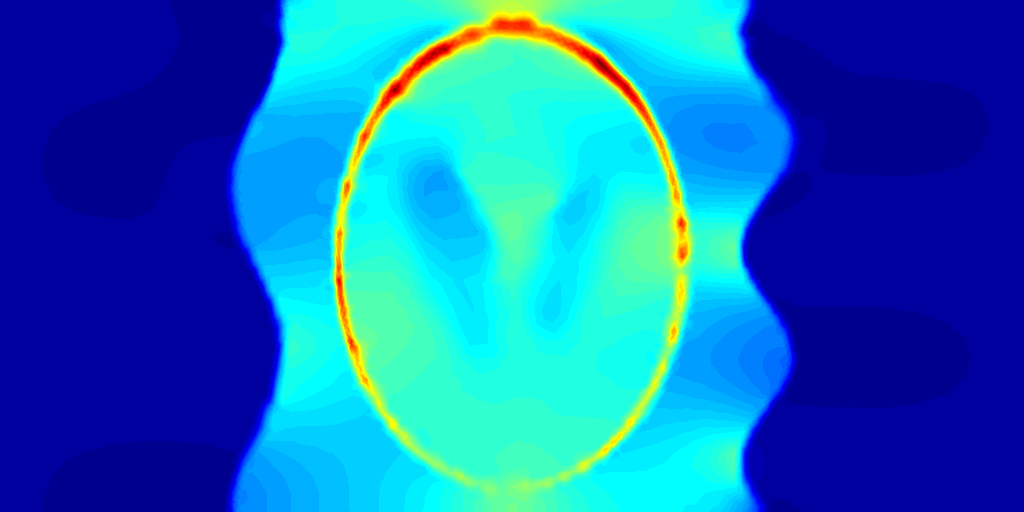};
\end{axis}
\end{tikzpicture}%
\caption{\label{min0} Reconstructed image without measurement noise.}
\end{center}
\end{figure}

\begin{figure}
\begin{center}
% This file was created by matlab2tikz v0.4.2.
% Copyright (c) 2008--2013, Nico Schlömer <nico.schloemer@gmail.com>
% All rights reserved.
% 
% The latest updates can be retrieved from
%   http://www.mathworks.com/matlabcentral/fileexchange/22022-matlab2tikz
% where you can also make suggestions and rate matlab2tikz.
% 
% 
% 
\begin{tikzpicture}

\begin{axis}[%
width=5.30555555555556in,
height=2.65277777777778in,
axis on top,
scale only axis,
xmin=0,
xmax=2,
ymin=0,
ymax=1,
title={Optimal control method},
colormap/jet,
colorbar,
point meta min=1.4236222413451,
point meta max=2.69934493292859
]
\addplot graphics [xmin=-0.000977517106549365,xmax=2.00097751710655,ymin=-0.000978473581213307,ymax=1.00097847358121] {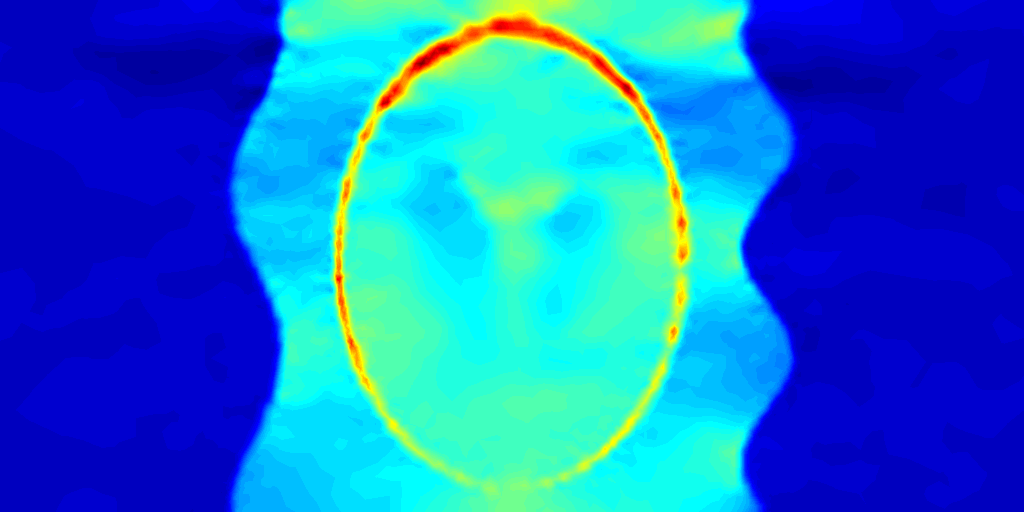};
\end{axis}
\end{tikzpicture}%
\caption{\label{min2} Reconstructed image with $2\%$ measurement noise.}
\end{center}
\end{figure}

\begin{figure}
\begin{center}
% This file was created by matlab2tikz v0.4.2.
% Copyright (c) 2008--2013, Nico Schlömer <nico.schloemer@gmail.com>
% All rights reserved.
% 
% The latest updates can be retrieved from
%   http://www.mathworks.com/matlabcentral/fileexchange/22022-matlab2tikz
% where you can also make suggestions and rate matlab2tikz.
% 
% 
% 
\begin{tikzpicture}

\begin{axis}[%
width=5.30555555555556in,
height=2.65277777777778in,
axis on top,
scale only axis,
xmin=0,
xmax=2,
ymin=0,
ymax=1,
title={Optimal control method},
colormap/jet,
colorbar,
point meta min=0.906385511876534,
point meta max=7.28404680791246
]
\addplot graphics [xmin=-0.000977517106549365,xmax=2.00097751710655,ymin=-0.000978473581213307,ymax=1.00097847358121] {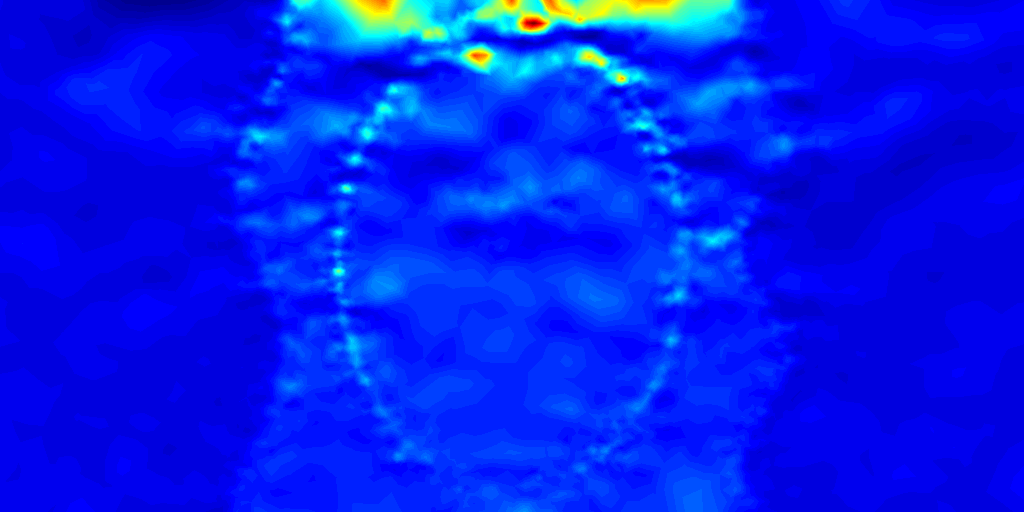};
\end{axis}
\end{tikzpicture}%
\caption{\label{min20} Reconstructed image with $20\%$ measurement noise.}
\end{center}
\end{figure}

The minimization procedure gives a decent qualitative reconstruction. The main interfaces are easy to see, yet this method, due to its regularizing effect, fails to show details in weaker contrasts zones. Figures \ref{min0}, \ref{min2}, and \ref{min20} show the reconstruction obtained with different measurement noise levels.

\subsubsection{The orthogonal field method}
\begin{figure}
\begin{center}
% This file was created by matlab2tikz v0.4.2.
% Copyright (c) 2008--2013, Nico Schlömer <nico.schloemer@gmail.com>
% All rights reserved.
% 
% The latest updates can be retrieved from
%   http://www.mathworks.com/matlabcentral/fileexchange/22022-matlab2tikz
% where you can also make suggestions and rate matlab2tikz.
% 
% 
% 
\begin{tikzpicture}

\begin{axis}[%
width=5.30555555555556in,
height=2.65277777777778in,
axis on top,
scale only axis,
xmin=0,
xmax=2,
ymin=0,
ymax=1,
title={Orthogonal field method},
colormap/jet,
colorbar,
point meta min=0.999367007051655,
point meta max=7.99500003130232
]
\addplot graphics [xmin=-0.000977517106549365,xmax=2.00097751710655,ymin=-0.000978473581213307,ymax=1.00097847358121] {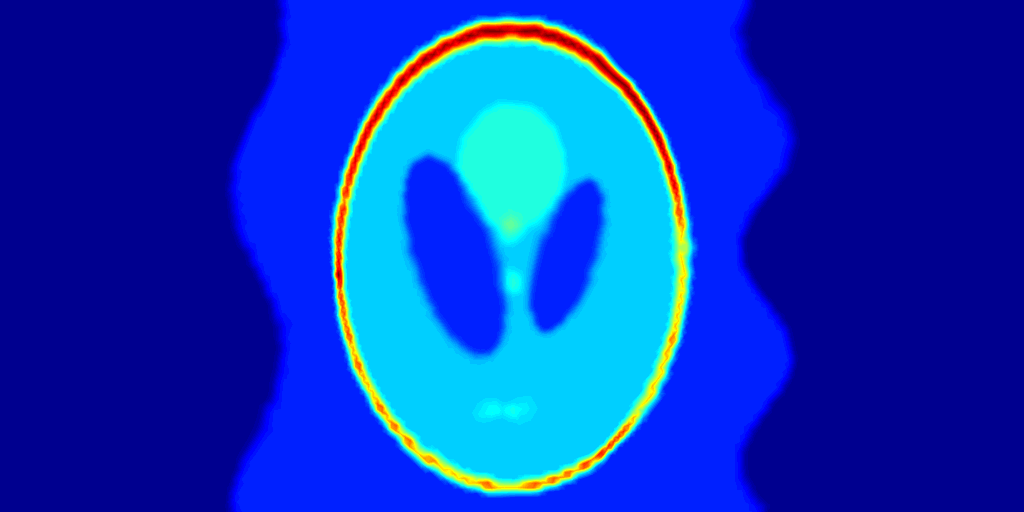};
\end{axis}
\end{tikzpicture}%
\caption{\label{visco0} Reconstructed image without measurement noise.}
\end{center}
\end{figure}

\begin{figure}
\begin{center}
% This file was created by matlab2tikz v0.4.2.
% Copyright (c) 2008--2013, Nico Schlömer <nico.schloemer@gmail.com>
% All rights reserved.
% 
% The latest updates can be retrieved from
%   http://www.mathworks.com/matlabcentral/fileexchange/22022-matlab2tikz
% where you can also make suggestions and rate matlab2tikz.
% 
% 
% 
\begin{tikzpicture}

\begin{axis}[%
width=5.30555555555556in,
height=2.65277777777778in,
axis on top,
scale only axis,
xmin=0,
xmax=2,
ymin=0,
ymax=1,
title={Orthogonal field method},
colormap/jet,
colorbar,
point meta min=0.844013668293187,
point meta max=8.33313009698197
]
\addplot graphics [xmin=-0.000977517106549365,xmax=2.00097751710655,ymin=-0.000978473581213307,ymax=1.00097847358121] {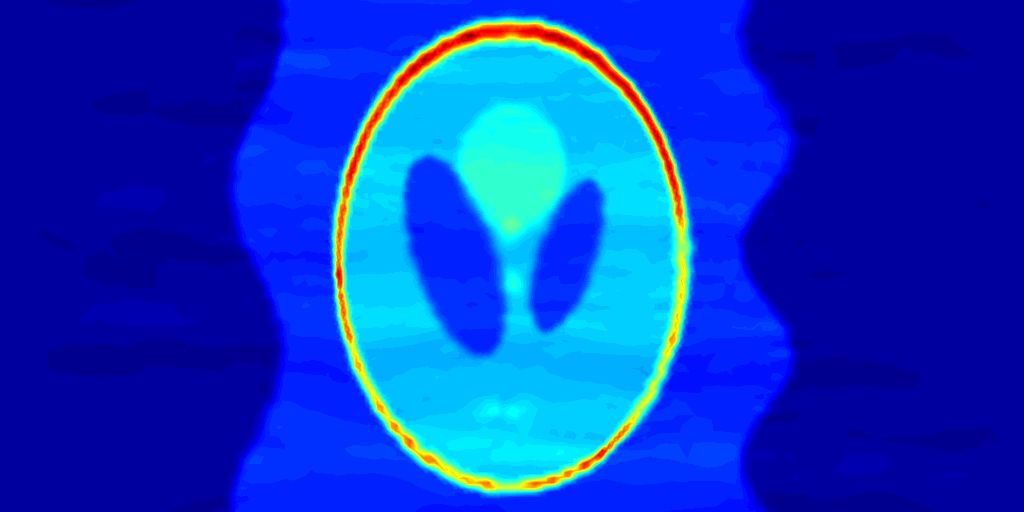};
\end{axis}
\end{tikzpicture}%
\caption{\label{visco2} Reconstructed image with $2\%$ measurement noise.}
\end{center}
\end{figure}

\begin{figure}
\begin{center}
% This file was created by matlab2tikz v0.4.2.
% Copyright (c) 2008--2013, Nico Schlömer <nico.schloemer@gmail.com>
% All rights reserved.
% 
% The latest updates can be retrieved from
%   http://www.mathworks.com/matlabcentral/fileexchange/22022-matlab2tikz
% where you can also make suggestions and rate matlab2tikz.
% 
% 
% 
\begin{tikzpicture}

\begin{axis}[%
width=5.30555555555556in,
height=2.65277777777778in,
axis on top,
scale only axis,
xmin=0,
xmax=2,
ymin=0,
ymax=1,
title={Orthogonal field method},
colormap/jet,
colorbar,
point meta min=0.498513511954329,
point meta max=13.1436521773483
]
\addplot graphics [xmin=-0.000977517106549365,xmax=2.00097751710655,ymin=-0.000978473581213307,ymax=1.00097847358121] {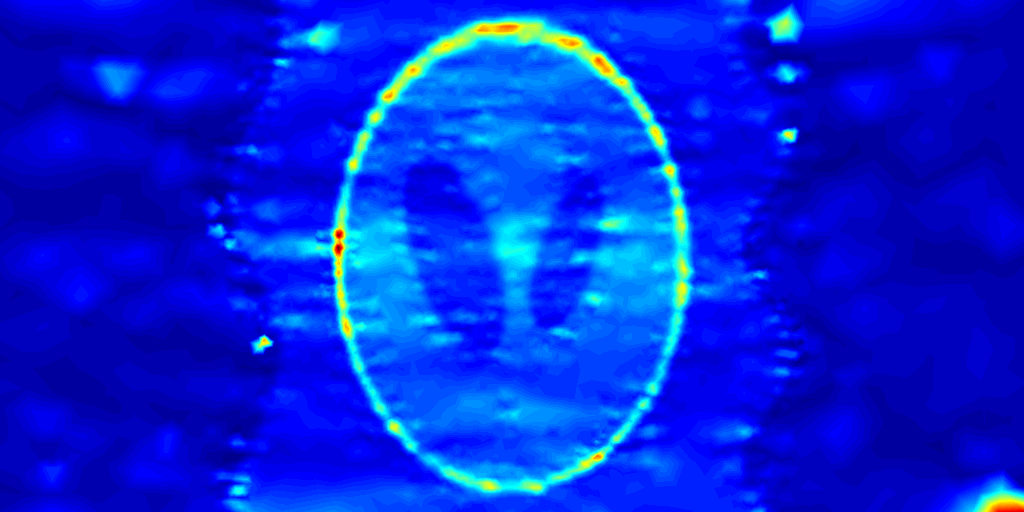};
\end{axis}
\end{tikzpicture}%
\caption{\label{visco20} Reconstructed image with $20\%$ measurement noise.}
\end{center}
\end{figure}

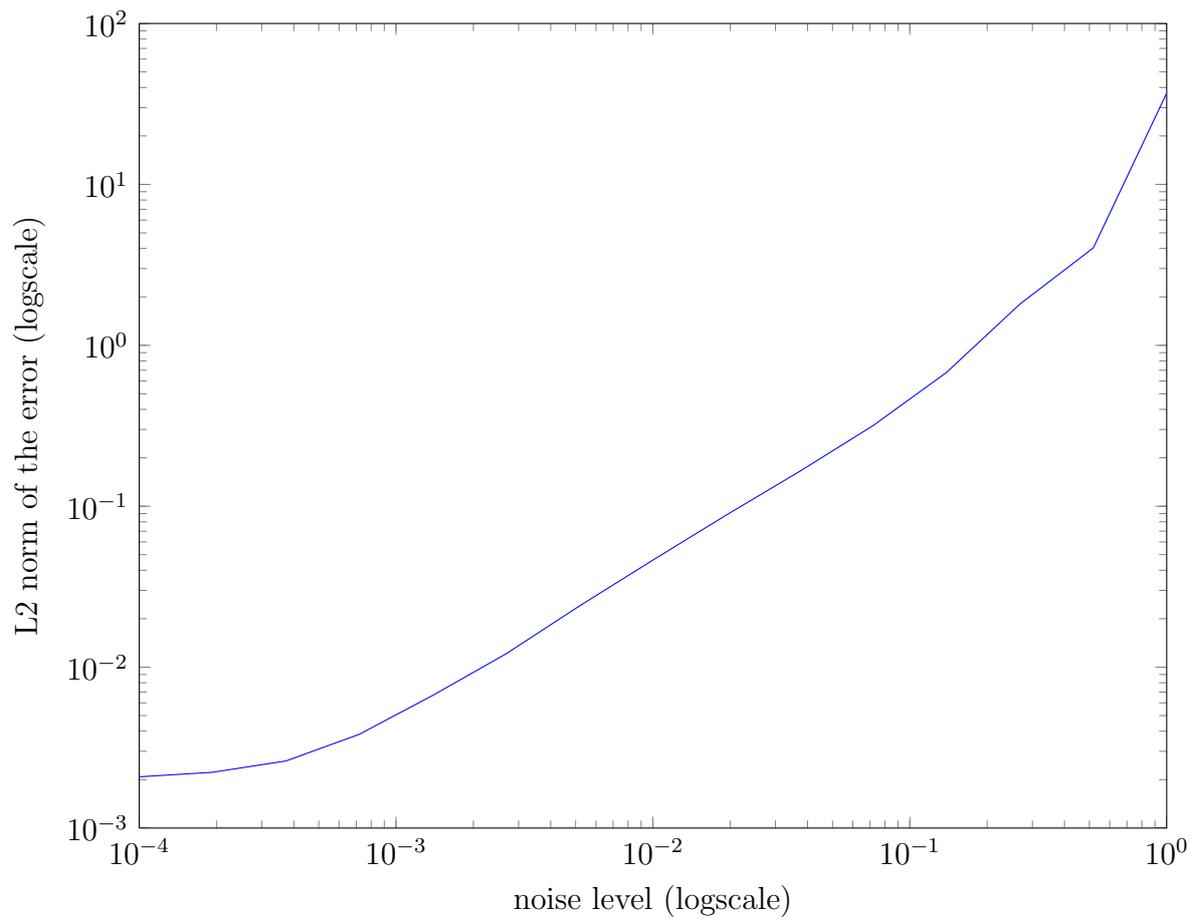
\begin{figure}
\begin{center}
% This file was created by matlab2tikz v0.4.2.
% Copyright (c) 2008--2013, Nico Schlömer <nico.schloemer@gmail.com>
% All rights reserved.
% 
% The latest updates can be retrieved from
%   http://www.mathworks.com/matlabcentral/fileexchange/22022-matlab2tikz
% where you can also make suggestions and rate matlab2tikz.
% 
% 
% 
\begin{tikzpicture}

\begin{axis}[%
width=5.323in,
height=4.20in,
scale only axis,
xmode=log,
xmin=0.0001,
xmax=1,
xminorticks=true,
xlabel={noise level (logscale)},
ymode=log,
ymin=0.001,
ymax=100,
yminorticks=true,
ylabel={L2 norm of the error (logscale)}
]
\addplot [
color=blue,
solid,
forget plot
]
table[row sep=crcr]{
0.0001 0.00208205294141109\\
0.000193069772888325 0.00222123329158285\\
0.000372759372031494 0.00261168517158364\\
0.000719685673001152 0.00382761408633793\\
0.00138949549437314 0.00667415162947117\\
0.00268269579527973 0.0121096075155522\\
0.00517947467923121 0.024015363811806\\
0.01 0.0462789966423204\\
0.0193069772888325 0.0881288601762636\\
0.0372759372031494 0.164735792657556\\
0.0719685673001152 0.317025262558384\\
0.138949549437314 0.67825957807824\\
0.268269579527972 1.80383889783423\\
0.517947467923121 4.02865355624427\\
1 36.9547629067761\\
};
\end{axis}
\end{tikzpicture}%
\caption{\label{noiselevel} $L^2$ norm of the error with respect to the noise level.}
\end{center}
\end{figure}

To find the solution of problem (\ref{eq_visco}), we fix $\varepsilon=10^{-3}$, and solve the equation on a uniform mesh on $\Omega$. We reconstruct an approximation of $\sigma$, and adapt the mesh to this first reconstruction. We do this procedure a couple of times in order to get refined mesh near the conductivity jumps. We can see that besides being computationally lighter than the minimization method, the orthogonal field method allows a quantitative reconstruction of $\sigma$ and shows details even in the low contrast zones.
It is relatively stable with respect to measurement noise. Figures \ref{visco0}, \ref{visco2}, and \ref{visco20} show the reconstruction with different measurement noise levels.
Figure \ref{noiselevel} shows the $L^2$ norm of the error with respect to measurement noise, with $\varepsilon$ fixed at $10^{-3}$. A smaller $\varepsilon$ increases the noise sensibility at higher noise levels, but also improves the details and reduces the smoothing effect of the $\varepsilon \Delta$ term in  (\ref{eq_visco}).

%\begin{figure}
%\begin{center}
%\input{figs/meshpng.tex}
%\caption{\label{imagemesh} Adapted mesh after $3$ iterations.}
%\end{center}
%\end{figure}


\section{Concluding remarks}

In this paper we have provided the mathematical basis of  ultrasonically-induced Lorentz force electrical impedance tomography. We have designed two efficient algorithms and tested them numerically. The resolution of the reconstructed images is fixed by the ultrasound wavelength and the width of the ultrasonic beam. The orthogonal field method performs much better than the optimization scheme in terms of both computational time and accuracy. 
 In a forthcoming work, we intend to generalize our approach for imaging anisotropic conductivities by ultrasonically-induced Lorentz force \cite{anisotropy}.   We will also propose an algorithm to find $\sigma \nabla U$ from the data function $\psi$ using (\ref{leading1})
and correct the leading-order approximation (\ref{leading2}). This will enhance the resolution of the reconstructed conductivity images. Another challenging problem under consideration is to interpret the high-frequency component of $M_{y,\xi}$ in terms of speckle conductivity contrasts.

 \end{document}